 \documentclass[reqno,12pt]{amsart}
 \usepackage{amsmath, latexsym, amsfonts, amssymb, amsthm, amscd, stmaryrd}

 \theoremstyle{plain}
 \newtheorem{theorem}{Theorem}%[section]
 \newtheorem*{theorem*}{Theorem}%[section]
 \newtheorem{lemma}{Lemma}
 \newtheorem{proposition}{Proposition}
 \newtheorem{corollary}{Corollary}
 \newtheorem{rem}{Remark}%[theorem]

 \newtheorem*{assumption*}{Assumption}

 \renewcommand{\tilde}{\widetilde}          % wider `tilde'
 \DeclareMathSymbol{\leqslant}{\mathalpha}{AMSa}{"36} % nicer `smaller or equal'
 \DeclareMathSymbol{\geqslant}{\mathalpha}{AMSa}{"3E} % nicer `larger or equal'
 \DeclareMathSymbol{\eset}{\mathalpha}{AMSb}{"3F}     % nicer `emptyset'
 \renewcommand{\leq}{\;\leqslant\;}                   % redef. of < or =
 \renewcommand{\geq}{\;\geqslant\;}                   % redef. of > or =
 \DeclareMathOperator*{\union}{\bigcup}       % \sum-like symbol for union
 \newcommand{\Z}{\mathbb{Z}}

 \title[Asymptotic direction for random walks in random environments]
 {Asymptotic direction for random walks in random environments}

 \author{Fran{\c c}ois Simenhaus}
 \thanks{Partially
    supported by CNRS (UMR 7599
 ``Probabilit{\'e}s et Mod{\`e}les
 Al{\'e}atoires'')}

\begin{document}

 \maketitle
 \begin{center}
 {\footnotesize \noindent
  Universit{\'e} Paris 7,\\
 Math{\'e}matiques, case 7012,\\ 2, place Jussieu, 75251 Paris, France}

 {\footnotesize \noindent e-mail: \texttt{simenhaus@math.jussieu.fr}}
 \end{center}

 \begin{abstract}
 In this paper we study the existence of an asymptotic direction for
 random walks in random i.i.d. environments (RWRE). We prove that if the
 set of directions
 where the walk is transient contains a  non empty open set,
 the walk admits an asymptotic direction.
 The main tool to obtain this result is the construction of 
 a renewal structure with cones. We also prove that RWRE admits
 at most two opposite asymptotic directions.
 \bigskip

 \noindent\textsc{Resum{\'e}:}
 Dans cet article, nous {\'e}tudions l'existence d'une direction
 asymptotique pour les marches al{\'e}atoires en milieu al{\'e}atoire
 i.i.d. (RWRE).
 Nous prouvons que si l'ensemble des directions dans lesquelles la marche est
 transiente contient un ouvert non vide, la marche admet une
 direction asymptotique. La construction d'une structure de
 renouvellement avec c{\^o}nes est le
 principal outil pour la preuve de ce r{\'e}sultat. Nous montrons aussi
 qu'une RWRE admet au plus 2 directions asymptotiques oppos{\'e}es.
 \end{abstract}

% %%%%%%%%%%%%%%%%%%%%%%%%%%%%%%%%%%%%%%%%%%%%%%%%%%%%%%%%%%%%%%%%%%%%%%%%%%%%%%%%%%%%
% %%%%%%%%%%%%%%%%%%%%%%%%%%%%%%%%%%%%%%%%%%%%%%%%%%%%%%%%%%%%%%%%%%%%%%%%%%%%%%%%%%%%
 \section{Introduction and results}
% %%%%%%%%%%%%%%%%%%%%%%%%%%%%%%%NOTAION%%%%%%%%%%%%%%%%%%%%%%%%%%%%%%%%%%%%%%%%%%%%%%
 In this paper, we give a characterization of random walks in
 random i.i.d environments having an asymptotic direction. We first
 describe the model that we will use. Fix a dimension
 $d\geq1$ (but think more particularly of the case where $d\geq 2$ because this work
 becomes obvious when $d=1$).
 Let $\mathcal{P}_+$ denote the $(2d-1)$ dimensional simplex,
 $\mathcal{P}_+=\{x\in[0,1]^{2d},\ \sum_{i=1}^{2d}x_i=1\}$.
 An environment $\omega$ in $\mathbb{Z}^d$ is an element
 of $\Omega:=\mathcal{P}_+^{\mathbb{Z}^d}$.
 For any environment $\omega$, $P_{x,\omega}$ denotes the Markov chain
 with state space $\mathbb{Z}^d$ and transition given by
 \begin{align*}
  P_{x,\omega}(X_0=x)&=1\qquad and\\
 P_{x,\omega}(X_{n+1}=z+e|X_n=z)&=\omega_z(e)\qquad(z\ \in \mathbb{Z}^d,\ e\ \in \mathbb{Z}^d\ s.t.\
 |e|=1,\ n\geq 0),
 \end{align*}
 where $|\cdot|$ denotes the Euclidean norm in $\mathbb{Z}^d$.\\
 For any law $\mu$ on $\mathcal{P}_+$, we define a random
 environment $\omega$ in $\mathbb{Z}^d$, random variable on $\Omega$ with law
 $\mathbb{P}:=\mu^{\otimes\mathbb{Z}^d}$.
 For any $x$ in $\mathbb{Z}^d$ and any fixed $\omega$, the
 law $P_{x,\omega}$ is called quenched
 law. The annealed law $P_x$ is defined on
 $\Omega\times(\mathbb{Z}^d)^{\mathbb{N}}$ by the
 semi-product $P_x:=\mathbb{P}\times P_{x,\omega}$. In this article,
 the law $\mu$ will verify the assumption of strict ellipticity
 \begin{equation*}
 \forall e \in \mathbb{Z}^d\
 s.t.\ |e|=1,\ \mathbb{P}-a.s.\quad\mu(\omega_0(e)>0)=1,
 \end{equation*}
 which is weaker than the usual uniform ellipticity (see Remark \ref{rem:ellipticity}).
 $S^{d-1}$ denotes the unit circle for the Euclidean norm.
 For any $\ell$ in $\mathbb{R}^d$, we define the set
 $A_{\ell}$ of
 transient trajectories in direction $\ell$
 \begin{equation*}
 A_{\ell}=\{\lim_{n\rightarrow +\infty} X_n\cdot\ell=+\infty\},
 \end{equation*}
 and for any $\nu$ in $S^{d-1}$, $B_{\nu}$ is defined as
 the set of trajectories having $\nu$ for asymptotic direction
 \begin{equation*}
 B_{\nu}=\{\lim_{n\rightarrow +\infty}\frac{X_n}{|X_n|}=\nu\}.
 \end{equation*}

 This model is well studied in the one dimensional case where many
 sharp properties of the walk are known. However in higher dimensions
 the behavior of the walk is much less
 well-understood. Particularly, the notion of asymptotic direction has
 been poorly studied. In this paper, we give a description of the class of walks
 having a unique asymptotic direction under the annealed measure
 (Theorem \ref{theorem1}).
 It means that the walk is
 transient and escapes to infinity in a direction which has a
 deterministic almost surely limit. We also prove that under the annealed
 measure, a RWRE admits at most two opposite asymptotic directions
 (Proposition \ref{proposition:vecteurs} and Corollary
 \ref{corollary:oppose}~).
 The proofs are based on renewal structure as in \cite{CZ_LLN} or \cite{SZ_LLN}.  

 The main difficulty to obtain an asymptotic direction for a transient
 walk is to control the fluctuations of the walk in the hyperplane
 transverse to transience direction.
 One way to control those fluctuations is to
 introduce the following assumption.
 \begin{assumption*}
 $\ell$ in $\mathbb{R}^d_*$ verifies assumption $(H)$
 if there exists a neighborhood $\mathcal{V}$ of $\ell$ such that
 \begin{equation*}
 \forall \ell' \in \mathcal{V},\quad P_0(A_{\ell'})=1.\qquad\qquad\textit{(H)}
 \end{equation*}
 \end{assumption*}
 When $(H)$ holds, we will note $\mathcal{V}$ the neighborhood
 given by the assumption.

 The main purpose of this article is to prove the following theorem.
 \begin{theorem}
 \label{theorem1}
 The following three statements are equivalent
 \begin{itemize}
 \item[\textit{i})]
 There exists a non empty open set  $\mathcal{O}$ of $\mathbb{R}^d$ such that
 \begin{equation*}
   \forall{\ell}\in\mathcal{O},\ \ \ \  P_{0}\left(A_{\ell}\right)=1.
 \end{equation*}
 \item[\textit{ii})]
  $\exists\nu\in S^{d-1}\ s.t.$
 \begin{equation*}
  P_{0}-a.s.,\quad \frac{X_n}{|X_n|}\xrightarrow[n\to\infty]{}\nu.
 \end{equation*}
  \item[\textit{iii})]
  $\exists\nu\in\mathbb{R}^d_*\ s.t.\ \forall\  \ell\in\mathbb{R}^d$
 \begin{equation*}
   \ell\cdot\nu>0\implies P_{0}\left(A_{\ell}\right)=1.
 \end{equation*}
 \end{itemize}
 \end{theorem}
 Using arguments similar to those applied in the proof of Theorem
 \ref{theorem1}, we also show

 \begin{proposition}
 \label{proposition:vecteurs}
 If $\nu$ and  $\nu'$ are two distinct vectors in $S^{d-1}$ such that
 $P_0(B_{\nu})P_0(B_{\nu'})>~0$ then $\nu'=-\nu$.
 \end{proposition}

 An obvious consequence of this proposition is the following corollary.

 \begin{corollary}
 \label{corollary:oppose}
 Under $P_0$, there are at most two asymptotic directions, in this
 case these two
 potential directions are opposite each other.
 \end{corollary}

 The class of walks admitting an asymptotic direction has been
 poorly studied so far. Theorem \ref{theorem1} gives a characterization of this
 class but also leaves unsolved some important problems related to this notion.
 First, we would like to compare the ballistic class with the class of
 walks admitting an asymptotic direction. As shown in Remark $(\ref{rem:lln})$,
 if a walk admits an asymptotic direction, it also satisfies a law of
 large numbers.
 However, the notion of asymptotic
 direction is of interest only in the non-ballistic case. Indeed, a non
 degenerate velocity contains more information (direction
 and speed) than the asymptotic direction (direction only) whereas
 in the non-ballistic case the asymptotic direction gives an interesting
 information of the behavior  which is not contained in the
 law of large numbers. It is known that in dimension $1$, the class of
 walks admitting an asymptotic direction (here it simply means
 transient) but a degenerate velocity is non empty. In higher dimension
 we have no example of such a walk and it might be possible that there
 is none.

 The class of ballistic walks has been the subject of many recent
 articles. In the first one \cite{SZ_LLN}, the authors provide a strong
 sufficient drift condition to obtain a ballistic law of large numbers,
 Kalikow's condition. Later, Sznitman improves on sufficient conditions in
 different works, \cite{S_TGAMA} is the more recent. He introduces in
 this paper the conditions $(T_{\gamma})\ (\gamma \in (0,1])$ that we
 will not recall, and the condition $(T')$
 defined as the realization of $(T_{\gamma})$ for any $\gamma \in
 (0,1)$. According to Corollary $5.3$ in \cite{S_NE}, this condition
 is strictly weaker than Kalikow's criterion for $d \geq 3$ and
 Sznitman gives an effective
 criterion to check it, that is the weakest condition known to assure a
 ballistic behavior. 
 It is also shown in $(1.13)$ of \cite{S_TGAMA} that $(T_{\gamma})$
 implies that the walk has an asymptotic direction
(and so using Theorem $1$ implies
 $(H)$), it is then natural to ask if $(H)$ is strictly weaker
 than $(T_{\gamma})$ or simply equivalent. This question is particularly
 interesting because $(T_{\gamma})$ is equivalent to $(T')$ for $\gamma
 \in (\frac{1}{2},1)$ when $d\geq 2$, and it is conjectured
 that they are equivalent for any
 $\gamma \in (0,1)$ (see \cite{S_TGAMA}, in particular Theorem $2.4$). An answer to this question
 could be a
 step toward comparing the ballistic class and the class of walks admitting an
 asymptotic direction and then would help us to solve the first
 problem above.

 A third problem is that we could not find
 a criterion to check assumption $(H)$, as we have for
 $(T_{\gamma})$. Such a criterion would be a great help to answer 
the two previous questions because condition $(H)$ is easy to
 understand but hard to verify. 

 % Finally, notice that Proposition \ref{proposition:vecteurs} and Corollary
%  \ref{corollary:oppose} are known to be true for velocities instead of
%  asymptotic directions (see \cite{BER}). However our
%  statements refine the former ones because if a walk
%  admits a non zero velocity $\nu$ then it also admits $\nu$ as
%  asymptotic direction. 

 Finally, notice that Corollary \ref{corollary:oppose} has a version
 for velocities instead of asymptotic directions. More precisely,
 $P_0$-almost surely, the limit of $\frac{X_n}{n}$
 belongs to a set $\{S_1,S_2\}$ such that
 if $S_1\neq0$ then $S_2=\lambda S_1$ for some
 $\lambda\leq 0$ (Theorem $1.1$ in \cite{BER}). However none of these
 two results can be deduced from the other one.
 
 The proofs of the results will be given in the second part of this
 paper. We finish this section with some notation which will be useful in the
 proofs. Denote by $\theta_n$ the time shift ($n$ natural number is
 the argument) and by $t_x$ the spatial shift ($x$ in $\mathbb{Z}^d$ is
 the argument). For any fixed $\ell$ in $\mathbb{R}^d_*$, we let $T_u$
 be the hitting time of the open half-space $\{x\in\mathbb{Z}^d,\ x\cdot
 \ell>u\}$
 \begin{equation*}
 T_u=\inf\{n>0,\ X_n \cdot \ell>u\},
 \end{equation*}
 and $D^\ell$ the return time of the walk behind the starting point
 \begin{equation*}
 D^\ell=\inf\{n> 0,\ X_n\cdot\ell\leq X_0 \cdot \ell\}.
 \end{equation*}
 Notice that these two definitions are quite different
 from those used in \cite{SZ_LLN}.\\
 We complete $\ell$ into an orthogonal
 basis $\left(e_2,\dots,e_d\right)$, such that
 for every $i$ in $\llbracket 2,d \rrbracket,\ |e_i|~=~1~.$\\
 For all $i\in \llbracket 2,d \rrbracket$ we define the following two
 vectors:
 \begin{equation}
\label{eq:defli}
 \ell'_i\left(\alpha\right)=\ell+\alpha e_i\quad \textrm{and}\quad \ell'_{-i}\left(\alpha\right)=\ell-\alpha e_i.
 \end{equation}
 For all positive real $\alpha$ we can define the convex
 cone $C\left(\alpha\right)$ by
 \begin{equation*}
 C\left(\alpha\right)=\bigcap_{i=2}^d\{x\in\mathbb{Z}^d,\ x\cdot
 \ell'_i\left(\alpha\right)\geq 0\quad \textrm{and}\quad x\cdot \ell'_{-i}\left(\alpha\right)\geq 0 \}.
 \end{equation*}
 We also define the exit time $D_{\alpha}^\ell$ of the cone
 $C\left(\alpha\right)$, shifted at the starting point of the walk,

 \begin{equation*}
 D_{\alpha}^\ell=\inf\{n\geq 0,\ \exists i \in \llbracket 2,d
 \rrbracket,\ X_n\cdot \ell'_i(\alpha)<X_0 \cdot \ell'_i(\alpha)\ \textrm{or}\ X_n\cdot \ell'_{-i}(\alpha)<X_0 \cdot \ell'_{-i}(\alpha)\}.
 \end{equation*}

 Notice that under $P_0$, $D_{\alpha}^\ell$ can also be defined in the
 following way,
 \begin{equation*}
 P_0-a.s.,\quad D_{\alpha}^\ell=\inf\{n\geq 0,\ X_n\notin C\left(\alpha\right)\}.
 \end{equation*}

% %%%%%%%%%%%%%%%%%%%%%%%%%%%%%%%%%%%%%%%%%%%%%%%%%%%%%%%%%%%%%%%%%%%%%%%%%%%%%%%%%%%%%%%%%%%%%%%%%%%%%%%%%%%%%%%%%%%%%%%%%%%%%%%%%%%%%%%%%%%%%%%%%%%%%%%%%%%%%%%%%%%%%%%%%%%%%%%%%%%%%%%%%%%%%%%%%%%%%%%%%%%%%%%%%%%%%%%%%%%%%%%%%%%%%%%%%%%%%%%%%%%%
 \section{Proofs}
% %%%%%%%%%%%%%%%%%%%%%%%%%%%%%%lemme1:cone positif%%%%%%%%%%%%%%%%%%%%%%%%%%%%%%%%%%%%%%%%%%%%%%%%%%%%%%%%%%%%%%%%%%%%%%%%%%%%%%%%%%%%%%%%%%%%%%%%%%%%%%%%%%%%%%
% \subsection{Lemma1}
 \begin{proof}[Proof of Theorem \ref{theorem1}]
 The first step of the proof is the following lemma, where it is proved
 that under $(H)$ the walk has a positive probability never to exit
 a cone $C(\alpha)$ for $\alpha$ small enough.
 \begin{lemma}
 \label{lemme:cone}
 Let $\ell$ be a vector in $\mathbb{R}^d$ satisfying $(H)$ then, for any choice of an  orthogonal basis
  $\left(\ell,e_2,\dots,e_d\right)$ with $|e_i|~=~1$
 for any $i$ in $\llbracket 2,d \rrbracket$,
 there exist some $\alpha_{0}>0$ such that,
 \begin{equation}
 \label{eq:preuve:cone}
   \forall{\alpha}\leq\alpha_{0}\ \ \ P_{0}\left(D_{\alpha}^{\ell}=\infty\right)>0.
 \end{equation}
 \end{lemma}
% %%%%%%%%%%%%%%%%%%%%%%%%%%%%%
% %%%%%%%%%%%%%%%%%%%%%%%%%%%%%
 \begin{proof}
 Fix a basis satisfying the assumption of the lemma.
 We will first show that there exists a random variable $\alpha_1>0$
 such that % $P_0(D_{\alpha_1}^\ell=\infty|D^{\ell}=\infty)=1$.
 \begin{equation}
 \label{eq:intermediaire}
 P_0(D_{\alpha_1}^\ell=\infty|D^{\ell}=\infty)=1.
 \end{equation}
 Since $\mathcal{V}$ is an open set, there exist
 some $\alpha_2>0$ such that for every
 $i\in\llbracket 2,d \rrbracket$:
 \begin{equation*}
  \ell'_i\left(\alpha_2\right)\in\mathcal{V}\ \textrm{and}\  \ell'_{-i}\left(\alpha_2\right)\in\mathcal{V}.
 \end{equation*}
 For these $\left(2d-2\right)$ directions, we use the renewal structure
 described in section $1$ of \cite{SZ_LLN}. The choice of the parameter
 $a$ in this structure has no importance and can be done arbitrarily.
 Remember that, for any fixed direction $\ell$, the first renewal time
 $\tau^{\ell}$ of \cite{SZ_LLN} is the first time the walk reaches a new
 record in direction $\ell$, and later never backtracks.
 \begin{rem}
 \label{rem:ellipticity}
 In \cite{SZ_LLN}, as in further references, uniform ellipticity is
 assumed. When we quote these articles, we have verified that this stronger
 assumption is not necessary or can be relaxed as in \cite{ZM}.
  \end{rem}
 Using $(H)$ we obtain that for each $i\in\llbracket 2,d \rrbracket$,
 \begin{equation*}
 P_0\left(A_{\ell'_i\left(\alpha_2\right)}\right)=P_0\left(A_{\ell'_{-i}\left(\alpha_2\right)}\right)=1.
 \end{equation*}
 From Proposition $1.2$ in \cite{SZ_LLN}:
 \begin{equation*}
 \tau^{\ell_2'\left(\alpha_2\right)}\vee\dots\vee\tau^{\ell_d'\left(\alpha_2\right)}\vee\tau^{\ell'_{-2}\left(\alpha_2\right)}\vee\dots\vee\tau^{\ell'_{-d}\left(\alpha_2\right)}<\infty\quad\quad
 P_0- a.s.
 \end{equation*}
 Using a proof very close to Proposition $1.2$ in \cite{SZ_LLN} (see
 also Theorem $3$ in \cite{KAL}) we obtain
 $P_0\left(D^{\ell}=\infty\right)~>~0$ and so:
 \begin{equation}
 \label{eq:tempsfini}
 \tau^{\ell_2'\left(\alpha_2\right)}\vee\dots\vee\tau^{\ell_d'\left(\alpha_2\right)}\vee\tau^{\ell_{-2}'\left(\alpha_2\right)}\vee\dots\vee\tau^{\ell'_{-d}\left(\alpha_2\right)}<\infty\quad\qquad
 P_0\left(\cdot |D^{\ell}=\infty\right)- a.s.
 \end{equation}
 We now define the following variables:
 \begin{eqnarray*}
 N&=&\inf\{n_0\geq 1,\ \forall n\geq n_0,\ X_n\in C\left(\alpha_2\right)\},\qquad(\inf\emptyset=+\infty)\\
 C&=&\inf_{1\leq n\leq N}X_n\cdot \ell,\\
 M&=&\sup_{1\leq n\leq N}\sum_{i=2}^d|X_n\cdot e_i|^2.
 \end{eqnarray*}
 From (\ref{eq:tempsfini}), it is clear that:
 \begin{eqnarray*}
 P_0\left(\cdot|D^{\ell}=\infty\right)-a.s.,&N&<\infty,\\
 &C>0\quad\textrm{and}&M<\infty.
 \end{eqnarray*}

 We now define $\alpha_1=\frac{C}{\sqrt{M}}\wedge\alpha_2$ (notice that $\alpha_1$ is
 random), using Cauchy-Schwarz
 inequality for $n\leq N$ and the definition of $N$ and $C(\alpha)$
 for $n\geq N$, we obtain:
 \begin{eqnarray*}
 P_0(\cdot|D^{\ell}=\infty)-a.s.,&\quad\forall i \in \llbracket 2,d
 \rrbracket\ \forall n\geq 0,\qquad&
 X_n\cdot \ell'_i\left(\alpha_1\right)=X_n\cdot
 \ell+\alpha_1\left(X_n\cdot e_i\right)\geq0,\\
 &and&X_n\cdot \ell'_{-i}\left(\alpha_1\right)=X_n\cdot
 \ell-\alpha_1\left(X_n\cdot e_i\right)\geq0.
 \end{eqnarray*}
 which ends the proof of (\ref{eq:intermediaire}).\\
 It is clear that
 \begin{equation*}
 \quad\alpha<\alpha'\quad\textrm{implies}\quad C\left(\alpha'\right)\subset C\left(\alpha\right),
 \end{equation*}
 and so\\
 \begin{equation*}
 \lim_{\alpha\rightarrow 0}P_0\left(\{D_{\alpha}^{\ell}=\infty\}\right)=P_0\left(\union_{\alpha>0}\{D_{\alpha}^{\ell}=\infty\}\right).
 \end{equation*}
 From (\ref{eq:intermediaire}), we have
 \begin{equation*}
 \union_{\alpha>0}\{D_{\alpha}^{\ell}=\infty\}\overset{P_0- a.s.}{=}\{D^{\ell}=\infty\}.
 \end{equation*}
 Since $P_0\left(D^{\ell}=\infty\right)>0$, this concludes the proof of
 Lemma \ref{lemme:cone}.
 \end{proof}
% %%%%%%%%%%%%%%%%%%%%%%%%%%%%%%%%%%%%%%%RENOUVELLEMENT%%%%%%%%%%%%%%%%%%%%%%%%%%%%%%%%%%%%%%%
% \subsection{Renewal structure with cones}
% \label{subsection:renouvellement}
 We will now construct a renewal structure in the same spirit as in
 \cite{SZ_LLN} or \cite{CZ_LLN}. The idea is to define a time where the walk reaches a new record in
 the direction $\ell$ and never goes out of a cone (also oriented in
 direction $\ell$) after. In \cite{SZ_LLN}, the walk moves from one slab to
 the next one, here, as in \cite{CZ_LLN} or \cite{CZ_FLU}, the walk will move from one cone to the next one.\\
 From Lemma \ref{lemme:cone}, we know that we can choose $\alpha$ small enough so that
 \begin{equation*}
 P_{0}\left(D_{\alpha}^{\ell}=\infty\right)>0.
 \end{equation*}
 We define now the two stopping time sequences $\left(S_k\right)_{k\geq 0}$and
 $\left(R_k\right)_{k\geq 0}$, and the sequence of successive maxima $\left(M_k\right)_{k\geq 0}$\\
 \begin{equation*}
  S_0=\inf\{n \geq 0,\ X_n \cdot \ell>X_0 \cdot \ell \},\quad
  R_0=D_{\alpha}^{\ell} \circ \theta_{S_0}+S_0,\quad M_0=\sup\{\ell\cdot X_n,\  0\leq n\leq R_0\}.
 \end{equation*}
 And for all $k\geq0$:
 \begin{equation*}
  S_{k+1}=T_{M_k},\quad
  R_{k+1}=D_{\alpha}^{\ell}\circ\theta_{S_{k+1}}+S_{k+1},\quad
  M_{k+1}=\sup\{\ell\cdot X_n,\ 0\leq n \leq R_{k+1}\},
 \end{equation*}
 \begin{equation*}
 K=\inf\{k\geq 0,\ S_k<\infty,\ R_k=\infty\}.
 \end{equation*}
 On the set ${K<\infty}$, we also define:
 \begin{equation*}
  \tau_1=S_K.
 \end{equation*}
 The random time $\tau_1$ is called the first cone renewal time, and
 will not be confused with $\tau^{\ell}$ introduced above.
 Under assumption $(H)$,
 \begin{equation}
 \label{ine}
 S_0\leq R_0 < S_1 \leq R_1< \dots<S_n\leq R_n < \cdots \leq\infty.
 \end{equation}
% %%%%%%%%%%%%%%%%%%%Proposition K fini%%%%%%%%%%%%%%%%%%%%%%%%%%%%%%%%%%
 \begin{proposition}Under assumption \textit{(H)},
 \label{propo:Kfini}
 \begin{equation*}
 P_0-a.s.\qquad K<\infty.
 \end{equation*}
 \end{proposition}
 \begin{proof}
 For all $k \geq 1$,
 \begin{eqnarray*}
 P_0(R_k<\infty)
 &=&\mathbb{E}[E_{0,\omega}[S_k<\infty,D_{\alpha}^{\ell}\circ\theta_{S_k}<\infty]]\\
 &=& \sum_{x\in
   \mathbb{Z}^d}\mathbb{E}[E_{0,\omega}[S_k<\infty,X_{S_k}=x,D_{\alpha}^{\ell}\circ\theta_{S_k}<\infty]].\\
 \end{eqnarray*}
 Using Markov property we obtain,
 \begin{equation*}
 P_0(R_k<\infty)=\sum_{x\in \mathbb{Z}^d}\mathbb{E}\left[E_{0,\omega}[S_k<\infty,X_{S_k}=x]E_{x,\omega}[D_{\alpha}^{\ell}<\infty]\right].
 \end{equation*}
 For every $x$ in $\mathbb{Z}^d$, the variables $E_{0,\omega}[S_k<\infty,X_{S_k}=x]$ and $E_{x,\omega}[D_{\alpha}^{\ell}<\infty]$ are
 respectively  $\sigma\{\omega_y(\cdot),\ \ell\cdot y< \ell\cdot x \}$ and $\sigma\{\omega_y(\cdot),\ y\in t_x\circ C\left(\alpha\right)\}$ measurable. As this
 two $\sigma$-fields are independent,
 \begin{eqnarray*}
 P_0(R_k<\infty)&=& \sum_{x\in \mathbb{Z}^d}E_0[S_k<\infty,X_{S_k}=x]E_x[D_{\alpha}^{\ell}<\infty]   \\
 &=&P_0(S_k<\infty)P_0(D_{\alpha}^{\ell}<\infty)\\
 &=&P_0(R_{k-1}<\infty)P_0(D_{\alpha}^{\ell}<\infty).\\
 \end{eqnarray*}
 By induction, we obtain,
 \begin{equation*}
 P_0(R_k<\infty)=P_0(D_{\alpha}^{\ell}<\infty)^{k+1}.
 \end{equation*}
 In view of (\ref{ine}), this concludes the proof.
 \end{proof}

 We now define a sequence of renewal time
 $\left(\tau_k\right)_{k\geq 1}$ by the following recursive relation:
 \begin{equation}
 \label{def:k+1}
 \tau_{k+1}=\tau_1\left(X.\right)+\tau_k\left(X_{\tau_1+.}-X_{\tau_1}\right).
 \end{equation}

 Using Proposition (\ref{propo:Kfini}), we have:
 \begin{equation*}
 \forall k\geq 0,\qquad \tau_k<\infty.
 \end{equation*}

 \begin{proposition}
 \label{prop:renouvellement}
 Under assumption $(H)$,\\
 $ \left( \left(X_{ \tau_1 \wedge\cdot} \right),\tau_1
   \right),\left(\left(X_{ \left( \tau_1+ \cdot \right) \wedge
   \tau_2}-X_{\tau_1}\right), \tau_2-\tau_1
   \right),\dots, \left( \left( X_{ \left( \tau_k + \cdot \right)
   \wedge \tau_{k+1}}-X_{\tau_k}\right),\tau_{k+1}-\tau_k
   \right)$
 are independent variables under $P_0$ and for $k\geq 1$,
 $\left(\left(X_{\left(\tau_k+\cdot\right)\wedge\tau_{k+1}}-X_{\tau_k}\right),\tau_{k+1}-\tau_k\right)$
 are distributed like
 $\left(\left(X_{\tau_1\wedge\cdot}\right),\tau_1\right)$ under
 $P_0\left(\cdot|D^\ell_{\alpha}=\infty\right)$.
 \end{proposition}

 The proof is similar to that of Corollary $1.5$ in \cite{SZ_LLN} and
 will not be repeated  here.\bigskip\\
% %%%%%%%%%%%%%%%%%%%%%%%%%%%%%%LEMME2:ESPERANCE X_tau%%%%%%%%%%%%%%%%%%%%%%%%%%%%%%%%%%%%
% \subsection{A renewal formula}
 For the classical renewal structure, Zerner proved that
 $E_0[X_{\tau_1} \cdot \ell]$ is finite and computes its value. We provide here the
 same result but for a renewal structure with cones.\\
 Fix a direction $\ell$ with integer coordinates $\left(a_1,\dots,a_d\right)$ such
 that their greatest common divisor,
 $\gcd\left(a_1,\dots,a_d\right)~=~1$.
 Assume that $(H)$ is satisfied for $\ell$.
 Complete $\ell$ in an orthogonal basis $\left(\ell,e_2,\dots,e_d\right)$ such that
 for every $i$ in $\llbracket 2,d \rrbracket,\ |e_i|~=~1~$. By Lemma
 \ref{lemme:cone}, we can choose $\alpha$ small enough so that
 $P_0(D_{\alpha}^{\ell}=\infty)>0$ and construct the associated renewal structure
 that is described above. %%in section \ref{subsection:renouvellement}.
 \begin{lemma}Under assumption \textit{(H)},
 \label{lemme:esperance}
 \begin{equation*}
 E_0[X_{\tau_1}\cdot\ell|D^{\ell}_{\alpha}=\infty]=\frac{1}{P_0\left(D_{\alpha}^{\ell}=\infty\right)}.
 \end{equation*}
 \end{lemma}
 \begin{proof}
 This proof follows an unpublished argument of M. Zerner but can be found in
 Lemma $3.2.5$ p$265$ of \cite{ZEI}. Since $\gcd\left(a_1,\dots,a_d\right)=1$, we have
 $\{x\cdot \ell,\ x\in\mathbb{Z}^d\}=\mathbb{Z}$. For all $i>0$,
 \begin{eqnarray}
 \label{ligne1}
 P_0\left(\{\exists k \geq 1,\ X_{\tau_k}\cdot\ell= i\}\right)  & = & \sum_{\{x\in\Z^d,x\cdot\ell=i\}}\mathbb{E}\left[E_{0,\omega}\left[\{X_{T_{i-1}}=x,D_{\alpha}^{\ell}\circ\theta_{T_{i-1}}=\infty\}\right]\right]
   \nonumber\\
 \label{ligne2}
 &=& \sum_{\{x\in\Z^d,x\cdot\ell=i\}}\mathbb{E}\left[E_{0,\omega}\left[X_{T_{i-1}}=x\right]E_{x,\omega}\left[D_{\alpha}^{\ell}=\infty\right]\right]\\
 \label{ligne3}
 &=&P_0\left(D_{\alpha}^{\ell}=\infty\right).
 \end{eqnarray}
 We used the strong Markov property in (\ref{ligne2}).\\
 In (\ref{ligne3}), we notice that $E_{0,\omega}\left(X_{T_{i-1}}=x\right)$ is $\sigma\{\omega_y(\cdot),\
 \ell\cdot y< \ell\cdot x \}$ measurable and $E_{x,\omega}\left(D_{\alpha}^{\ell}=\infty\right)$ is
 $\sigma\{\omega_y(\cdot),\ y\in t_x\circ C\left(\alpha\right)\}$ measurable and that those
 two $\sigma$-fields are independent. We will now compute the same value in another way.
 \begin{align*}
 \lim_{i\rightarrow\infty}P_0(\{\exists k &  \geq 1,\ X_{\tau_k}\cdot
   \ell= i\})   =  \lim_{i\rightarrow\infty}P_0\left(\{\exists k
   \geq 2,\ X_{\tau_k}\cdot \ell= i\}\right)\\
 &=\lim_{i\rightarrow\infty}\sum_{n\geq 1}P_0\left(\{\exists k \geq
   2,\ \left(X_{\tau_k}-X_{\tau_1}\right)\cdot \ell= i-n,X_{\tau_1}\cdot \ell=n\}\right)\\
 &=\lim_{i\rightarrow\infty}\sum_{n\geq 1}P_0\left(\{\exists k \geq
   2,\ \left(X_{\tau_k}-X_{\tau_1}\right)\cdot \ell=
   i-n\}\right)P_0\left(X_{\tau_1}\cdot \ell=n\right).
 \end{align*}
 Notice also that the first equality is true because $P_0\left(X_{\tau_1}\cdot\ell>i\right)\rightarrow 0
 \left(i\rightarrow\infty\right)$. Using now the renewal theorem (Corollary $10.2$ p$76$ in \cite{THO}) we obtain
 \begin{equation*}
 \lim_{i\rightarrow\infty}P_0\left(\{\exists k \geq 2,\
   \left(X_{\tau_k}-X_{\tau_1}\right)\cdot \ell=
   i-n\}\right)=\frac{1}{E_0[\left(X_{\tau_2}-X_{\tau_1}\right)\cdot \ell]}.
 \end{equation*}
 The dominated convergence theorem leads to
 \begin{equation*}
 \lim_{i\rightarrow\infty}P_0\left(\{\exists k \geq 1,\ X_{\tau_k}\cdot\ell=
   i\}\right)=\frac{1}{E_0[\left(X_{\tau_2}-X_{\tau_1}\right)\cdot \ell]}.
 \end{equation*}
 Comparing this result with (\ref{ligne3}), we easily obtain Lemma \ref{lemme:esperance}.
 \end{proof}

% %%%%%%%%%%%%%%%%%%%%%%%%%%%%%%%%%%%%FIN THEOREME%%%%%%%%%%%%%%%%%%%%%%%%%%%%%%%%%%%%%%%%%%%%%%%%
% \subsection{Proof of theorem \ref{theorem1}}
 We have now all the tools to prove Theorem \ref{theorem1}.
 We will first use the two lemmas to prove that $\textit{i})$ implies $\textit{ii})$.\\
 We choose $\ell$ with rational coordinates in the open set
 $\mathcal{O}$. It is
 clear that $\ell$ satisfies assumption $(H)$. Actually, we can also assume,
 without loss of generality, that $\ell$ has integer
 coordinates and that their greatest  common divisor is $1$.
 Indeed, there is $\lambda$ rational such that $\lambda\ell$ has integer
 coordinates with greatest common divisor equal to $1$, and of
 course, $\lambda\ell$ also satisfies $(H)$.\\

 We complete $\ell$ into an orthogonal
 basis $\left(e_2,\dots,e_d\right)$, such that
 for every $i$ in $\llbracket 2,d \rrbracket,\ |e_i|~=~1~.$
  Using Lemma \ref{lemme:cone}, we choose $\alpha$ small
 enough so that $P_0\left(D_{\alpha}^{\ell}=\infty\right)~>~0$.
 We can now use the renewal structure with cones %described in
% section (\ref{subsection:renouvellement})
 and we have from Lemma \ref{lemme:esperance} that:
 \begin{equation*}
 E_0\left[X_{\tau_1}\cdot \ell|D_{\alpha}^\ell=\infty\right]=\frac{1}{P_0\left(D^\ell_{\alpha}=\infty\right)}<\infty.
 \end{equation*}
 From the definition of the cone renewal structure, there is some
 constant $c(\alpha)>0$ such that,
 $P_0(\cdot|D_{\alpha}^{\ell}=\infty)-a.s.$, for any time $n$,
 \begin{equation}
 \label{eq:controle_norme}
 |X_n|\leq c(\alpha)X_n\cdot \ell,
 \end{equation}
% \begin{equation}
% \label{eq:controle_norme}
% \exists\  c(\alpha) \in\mathbb{R^+} s.t.\quad
% P_0(\cdot|D_{\alpha}^{\ell}=\infty)-a.s., \forall n \geq 0,
% \quad |X_n|\leq c(\alpha)\
% X_n\cdot \ell,
% \end{equation}
 and so using Lemma \ref{lemme:esperance},
 \begin{equation}
 \label{esperance_finie}
 E_0\left[|X_{\tau_1}|\vert D_{\alpha}^\ell=\infty\right]<\infty.
 \end{equation}
 We can now apply the law of large numbers, and obtain
 \begin{equation}
\label{eq:llninter}
 \frac{X_{\tau_k}}{k}\xrightarrow[n\to\infty]{}
E_0\left[X_{\tau_1}|D^{\ell}_{\alpha}=\infty\right]\quad P_0-a.s.
 \end{equation}
 As $|E_0\left[X_{\tau_1}|D_{\alpha}^\ell=\infty\right]|>0$,
 \begin{equation}
 \label{eq:lln}
 \frac{X_{\tau_k}}{|X_{\tau_k}|}\xrightarrow[n\to\infty]{}
\frac{E_0\left[X_{\tau_1}|D_{\alpha}^{\ell}=\infty\right]}
{|E_0\left[X_{\tau_1}|D_{\alpha}^\ell=\infty\right]|}\stackrel{(def)}{=}\nu\quad
P_0-a.s.
 \end{equation}
 To complete the proof, we have to control the behavior of the walk between the renewal
 times. For each natural $n$, we introduce the index $k\left(n\right)$ such that,
 \begin{equation*}
 \label{eq:def_tau}
 \tau_{k\left(n\right)}\leq n<\tau_{k\left(n\right)+1}.
 \end{equation*}
 Recall that if $\left(Z_n\right)$ is an i.i.d. sequence of variables with finite
 expectation, Borel-Cantelli Lemma assures that $\frac{Z_n}{n}$ converges
 almost surely to $0$. From (\ref{eq:controle_norme}), Lemma~\ref{lemme:esperance}
 and Proposition \ref{prop:renouvellement},
 the sequence
 $\displaystyle(\sup_n|X_{\tau_{k}+n\wedge\tau_{k+1}}-X_{\tau_k}|)_{k\geq
 1}$
 is i.i.d with finite expectation and
 we can apply the previous remark to obtain:
 \begin{equation}
 \label{eq:controle}
 \frac{\displaystyle\sup_n|X_{\tau_{k}+n\wedge\tau_{k+1}}-X_{\tau_k}|}{k}
\xrightarrow[k\to\infty]{}0\qquad P_0-a.s.
 \end{equation}
 Using equation (\ref{eq:lln}) and (\ref{eq:controle}), we study the convergence
 in $ii)$,
 \begin{equation}
 \label{eq:decompo}
 \frac{X_n}{|X_n|}=\frac{X_n-X_{\tau_{k\left(n\right)}}}{|X_n|}+\frac{X_{\tau_{k\left(n\right)}}}{k\left(n\right)}\frac{k\left(n\right)}{|X_n|}.
 \end{equation}
 By Proposition \ref{propo:Kfini} and (\ref{def:k+1}),
 \begin{equation*}
 k\left(n\right)\xrightarrow[n\to\infty]{}\infty\qquad P_0-a.s.
 \end{equation*}
 As $|X_n|\geq k \left( n \right)$, (\ref{eq:controle}) leads to:
 \begin{equation}
 \label{eq:intermediaire2}
 \frac{X_n-X_{\tau_{k\left(n\right)}}}{|X_n|}\xrightarrow[n\to\infty]{}0\qquad P_0-a.s.
 \end{equation}
 To control the second term in (\ref{eq:decompo}), we simply write
 \begin{equation*}
 \frac{|X_{\tau_{k\left(n\right)}}|}{k\left(n\right)}- \frac{|X_n-X_{\tau_{k\left(n\right)}}|}{k\left(n\right)} \leq \frac{|X_n|}{k\left(n\right)} \leq \frac{|X_{\tau_{k\left(n\right)}}|}{k\left(n\right)}+\frac{|X_n-X_{\tau_{k\left(n\right)}}|}{k\left(n\right)}.
 \end{equation*}
 Using (\ref{eq:llninter}) and (\ref{eq:controle}), we obtain,
 \begin{equation*}
 \frac{|X_n|}{k\left(n\right)}\xrightarrow[n\to\infty]{}|E_0\left[X_{\tau_1}\cdot \ell|D_{\alpha}^\ell=\infty\right]|\qquad P_0-a.s.
 \end{equation*}
 We finally obtain the desired convergence :
 \begin{equation*}
 \frac{X_n}{|X_n|}\xrightarrow[n\to\infty]{}\nu=\frac{E_0\left[X_{\tau_1}|D_{\alpha}^\ell=\infty\right]}{|E_0\left[X_{\tau_1}|D_{\alpha}^\ell=\infty\right]|}\qquad P_0-a.s.
 \end{equation*}
% %%%%%%%%%%%%%%%%%%%%%%%%%%%%%%%%Le reste du theor{\`e}me%%%%%%%%%%%%%%%%%%%%%%%%%%%%%%%%%%%%%%%%%%%%
 The end of the proof of Theorem \ref{theorem1} is easy: it is obvious that $\textit{iii})$ implies $\textit{i})$ and so we just have to show
 that $\textit{ii})$ implies $\textit{iii})$.\\
 Let $\ell$ be a direction such that $\ell \cdot \nu>0$. It is known since
 \cite{SZ_LLN} (Lemma $1.1$) that $P_0\left(A_\ell\cup A_{-\ell}\right)$ follows a
 $0-1$ law under assumption of uniform ellipticity, but we use
 here Proposition $3$ in \cite{ZM} where the same result is proved under the weaker assumption of
 strict ellipticity.\\
 If $P_0\left(A_\ell\cup A_{-\ell}\right)=0$, it is known that the walk under $P_0$
 oscillates,
 \begin{equation*}
 \limsup_{n\rightarrow\infty}X_n\cdot\ell=-\liminf_{n\rightarrow\infty}X_n\cdot\ell=+\infty\qquad P_0-a.s.
 \end{equation*}
 This is not possible in view of $\textit{ii)}$ and so $P_0\left(A_\ell\cup A_{-\ell}\right)=1$.\\
 But because of $\textit{ii})$, $P_0\left(A_{-\ell}\right)=0$, and we can conclude
 \begin{equation*}
 P_0\left(A_\ell\right)=1.
 \end{equation*}
 \end{proof}
 \begin{rem}
 \label{rem:lln}
 From the proof of Theorem \ref{theorem1}, we know that if a walk has
 an asymptotic direction, we can construct a renewal structure with
 cones and $E[X_{\tau_1}|D_{\alpha}^{\ell}=\infty]$ is finite. We can then easily derive a law
 of large numbers, namely
 \begin{equation*}
 \frac{X_n}{n}\xrightarrow[n\to\infty]{} \frac{E_0[X_{\tau_1}|D_{\alpha}^{\ell}=\infty]}
 {E_0[\tau_1|D_{\alpha}^{\ell}=\infty]}\stackrel{(def)}{=}\mu\qquad P_0-a.s.
 \end{equation*}
 However this limit can be zero (if and only if
 $E_0[\tau_1|D_{\alpha}^{\ell}=\infty]=+\infty$ ) 
 and, in this case, the asymptotic
 direction is an interesting information about the walk's behaviour.
 \end{rem}
 \begin{rem}
 \label{rem:classe}
If a walk admits an asymptotic direction, we know that the walk is
transient in any direction $\ell$ satisfying $\ell\cdot\nu>0$.
For any such $\ell$, we can consider a slab renewal
structure defined as the cone renewal structure except that in the
definitions of $(R_k)_{k\geq 0},(S_k)_{k\geq 0}$ and $(M_k)_{k\geq 0}$,
$D^\ell_\alpha$ is replaced by
$\tilde{D}^\ell=\inf\{n>0,X_n\cdot\ell<X_0\cdot\ell\}$ (this construction
is very similar to that of \cite{SZ_LLN}). For any $k\geq 1$,
$\tau^\ell_k$ will denote the $k$-th slab renewal time.
The variables $(X_{\tau^\ell_{k+1}}-X_{\tau^\ell_k})_{k\geq 1}$ are
i.i.d., and the purpose of this remark is to show that
the expectation of their norm,
$E_0[|X_{\tau_1^{\ell}}||\tilde{D}^{\ell}=\infty]$ is finite.
From Corollary $3$ in \cite{KES_LLN}, it is enough to prove
that $(\frac{X_{\tau^\ell_k}}{k})_{k\geq 1}$ is bounded $P_0-$almost
surely\footnote{The result in \cite{KES_LLN} is for $d=1$. Considering all the
  coordinates separately, the claim in arbitrary dimension follows.}.
 For any
$k\geq 1$, we introduce
\begin{equation*}
J(k)=\sup\{j\geq 1,\tau_j \leq \tau^\ell_k\}\qquad (\sup\emptyset=0).
\end{equation*}
As $P_0$-almost surely $(\tau_k^\ell)_{k\geq 1}$
is increasing, we have
\begin{equation}
\label{Jlim}
\lim_{k\rightarrow\infty}J(k)=\infty\qquad P_0-a.s.
\end{equation}
Notice that a cone renewal time is also a slab renewal
time and as a consequence
\begin{equation}
\label{Jinf}
 J(k)\leq k\qquad P_0-a.s.
\end{equation}
$P_0$-a.s. for large $k$ so that $J(k)\geq 1$:
\begin{equation*}
\frac{X_{\tau^\ell_k}}{k}=
\frac{X_{\tau^\ell_k}-X_{\tau_{J(k)}}}{k}+\frac{X_{\tau_{J(k)}}}{k}.
\end{equation*} 
The norm of the first term is bounded
by $\frac{c(\alpha)|X_{\tau_{J(k)+1}}\cdot \ell - X_{\tau_{J(k)}}\cdot
  \ell|}{J(k)}$. From the same argument as in (\ref{eq:controle}), $\frac{c(\alpha)|X_{\tau_{j+1}}\cdot \ell - X_{\tau_j}\cdot
  \ell|}{j}$ converges $P_0$-almost surely to $0$, and using (\ref{Jlim}), we obtain
the $P_0$-almost surely convergence of the first term to $0$.
Rewriting the second term 
\begin{equation*}
\frac{X_{\tau_{J(k)}}}{k}=\frac{X_{\tau_{J(k)}}}{J(k)}\frac{J(k)}{k},
\end{equation*}
from (\ref {Jinf}),(\ref{Jlim}) and (\ref{eq:llninter}) we obtain that the second term
is almost surely bounded as $k$ goes to infinity.

% We can deduce that the i.i.d. variables $(X_{\tau^\ell_{k+1}}-X_{\tau^\ell_k})_{k\geq 1}$
% satisfy a law of large number with finite limit
% and then the expectation of their norm, $E_0[|X_{\tau_1^{\ell}}||D^{\ell}=\infty]$ is finite.
% As a consequence if a walk admits an asymptotic direction, for any
% $\ell$ such that $\nu\cdot\ell>0$, the walk is transient in the
% direction $\ell$ and  $E_0[|X_{\tau_1^{\ell}}||D^{\ell}=\infty]<\infty$.

% $X_{\tau_1}$ has the same law under $P_0(\cdot|D^\ell=\infty)$
% that $\sum_{k=0}^J Y_k$ where $N$ is a geometrical variable with
% success parameter $\frac{P_0(D_{\alpha}^\ell=~\infty~)~}{P_0(D^\ell=\infty)}$ 
% independent from $(Y_k)_{k\geq 0}$ which is an i.i.d. family of
% variable with same law as
% $X_{\tau^{\ell}}$
% under $P_0(\cdot|\{D^\ell=\infty\}\cap\{D^\ell_{\alpha}<\infty\})$.
% As consequence if a walk admits an asymptotic direction
% then $E[|X_{\tau^{\ell}}||D^{\ell}=\infty]<\infty$.

 % The class of walks admitting an asymptotic direction is exactly the class of
%  transient walks (in a direction $\ell$) such that
%  $E[|X_{\tau^{\ell}}||D^{\ell}=\infty]<\infty$. 
%  For one of the inclusion just notice that $\tau_1$
%  is also one of the hyperplane renewal time. For the other one,
%  use the end of the proof of \textit{i)} implies \textit{ii)} in
%  Theorem \ref{theorem1} (from (\ref{esperance_finie})) with $\tau^{\ell}$
%  instead of $\tau_1$.  
 \end{rem}
% \begin{rem}
% It is not easy to check \textit{i)} (or \textit{(H)}). However we
% can find in previous papers some examples of walks admitting an
% asymptotic direction and then satisfying \textit{i}). Of course all
% criteria assuring a ballistic law of large numbers work like Kalikow's
% condition in \cite{KAL} or more recently condition $(T')$ of Sznitman in
% \cite{S_TGAMA}. The condition $(T)_{\gamma}$ described in \cite{S_TGAMA}
% also works (see $(1.13)$ in \cite{S_TGAMA}) and it is not proved that it
% implies a law of large numbers. More generally, the previous remark
% gives criteria to describe exhaustively the class of
% walks admitting an asymptotic direction.
% \end{rem}

%%%%%%%%%%%%%%%%%%%%%%%%%%%%%%%%%%%%%%%%%%%%%%%%%%%%%%%%%%%%%%%%%%%%%%%%%%%%%%%%%%%%%%%%%%%%
%%%%%%%%%%%%%%%%%%%%%%%%%%%%%%%%%Preuve proposition 1%%%%%%%%%%%%%%%%%%%%%%%%%%%%%%%%%%%%%%%
%%%%%%%%%%%%%%%%%%%%%%%%%%%%%%%%%%%%%%%%%%%%%%%%%%%%%%%%%%%%%%%%%%%%%%%%%%%%%%%%%%%%%%%%%%%%

\begin{proof}[Proof of Proposition \ref{proposition:vecteurs}]

Suppose that the proposition is false and call $\nu$ and $\nu'$ two
vectors of  $S^{d-1}$ different and non opposite such that
$P_0(B_{\nu})P_0(B_{\nu'})>~0$, then we will show that
\begin{equation}
\label{but}
\exists \nu_0\ \textrm{such that}\ P_0(B_{\nu_0}|B_{\nu}\cup B_{\nu'})=1,
\end{equation}
what establishes, of course, a contradiction.\\
 \\
First, notice that for $\nu\in S^{d-1}$ with $P_0(B_{\nu})>0$ we have,
\begin{equation}
\label{eq:cond}
\forall \ell\in \mathbb{R}^{d}\ \textrm{such that}\ \ell\cdot\nu>0,\quad P_0(A_{\ell}|B_{\nu})=1.
\end{equation}
Indeed, from the $0-1$ law $P_0(A_{\ell}\cup
A_{-\ell})=1$ (just notice that the walk does not oscillate along
direction $\ell$ on $B_{\nu}$, event of positive probability), but since
$B_{\nu}\subset\{\exists N\ s.t.\ \forall n\geq N,\
X_n\cdot\ell>0\}$, we have $B_{\nu}~\subset~A_{\ell}$, $P_0$-almost surely,
which implies (\ref{eq:cond}).

The set $\mathcal{V}=\{\ell\in\mathbb{R}^{d},\
\ell\cdot\nu>0\}\cap\{\ell\in\mathbb{R}^{d},\ \ell\cdot\nu'>0\}$ is
non empty and open and from (\ref{eq:cond}) it has the property,
% \begin{equation}
% \label{newh}
% \forall \ell\in E,\qquad P_0(A_{\ell}|B_{\nu}\cup B_{\nu'})=1.
% \end{equation}
% Fix a vector $\ell_0$ in $E$, it is obvious that there exists a
% neighborhood $\mathcal{V}$ of $\ell_0$ such that,
\begin{equation*}
\forall \ell \in \mathcal{V},\qquad P_0(A_{\ell}|B_{\nu}\cup B_{\nu'})=1.\qquad\qquad\textit{(H')}
\end{equation*}
From now on, we fix $\ell_0$ in $\mathcal{V}$.
The property $\textit{(H')}$ is similar to assumption $(H)$ and the
proof of (\ref{but}) will be adapted from that of Theorem 
\ref{theorem1}. In fact, we will use the cone renewal structure on
$A_{\ell_0}$ and show
\begin{equation}
\label{but2}
\exists \nu_0 \ s.t.\ P_0(B_{\nu_0}|A_{\ell_0})=1,
\end{equation}
which is stronger that (\ref{but}).
 Before the proof, notice that if $P_0(B_{\nu}\cup
B_{\nu'})=1$, we can easily conclude using Theorem~\ref{theorem1},
but the proof is not that obvious if
${P_0(B_{\nu}\cup B_{\nu'})<1}$.\\

In order to construct a renewal structure with cones on the event
$A_{\ell_0}$ of positive probability, we have to show that there
exists $\alpha>0$ such that $P_0(D^{\ell_0}_{\alpha}=\infty)>0$. We will use a
proof very close to the one of Lemma~\ref{lemme:cone} except that we will work
on the event $\{B_{\nu}\cup B_{\nu'}\}$ and show the stronger
result $P_0(D^{\ell_0}_{\alpha}=\infty,B_{\nu}\cup B_{\nu'})>0$.
Our first step will be to prove that
\begin{equation}
\label{positif}
P_0(D^{\ell_0}=\infty,B_{\nu}\cup B_{\nu'})>0.
\end{equation}
We argue by contradiction and assume the left hand side of
 (\ref{positif}) is zero. By translation-invariance,
 it holds $P_x(D^{\ell_0}=\infty,B_{\nu}\cup B_{\nu'})=0$ for any $x$ in
$\mathbb{Z}^d$, which means
\begin{equation*}
\forall x \in \mathbb{Z}^d,\ \mathbb{P}-a.s.,\quad P_{x,\omega}(D^{\ell_0}=\infty,B_{\nu}\cup B_{\nu'})=0,
\end{equation*}
and also
\begin{equation*}
\mathbb{P}-a.s.,\ \forall x \in \mathbb{Z}^d,\quad P_{x,\omega}(D^{\ell_0}=\infty,B_{\nu}\cup B_{\nu'})=0.
\end{equation*}
We denote by $(D^{\ell_0})^n$ the $n$-th backtrack time of the walk, defined
by the following recursive relation
\begin{eqnarray*}
(D^{\ell_0})^1&=&D^{\ell_0},\\
(D^{\ell_0})^n&=&D^{\ell_0}\circ
\theta_{(D^{\ell_0})^{n-1}}1_{\{(D^{\ell_0})^{n-1}<\infty\}}+\infty 1_{\{D^{\ell_0})^{n-1}=\infty\}},\quad \forall n
\geq 2.
\end{eqnarray*}
Notice that for any $n\geq 1$, 
\begin{align*}
\{&(D^{\ell_0})^{n}<\infty,(D^{\ell_0})^{n+1}=\infty,B_{\nu}\cup B_{\nu'}\}
=\\
&\bigcup_{x\in\mathbb{Z}^d}\bigcup_{m\geq
  1}\{(D^{\ell_0})^n<\infty,(D^{\ell_0})^n=m,X_m=x\}
\bigcap\{D^{\ell_0}\circ\theta_m=\infty,\lim_{k\rightarrow\infty} \frac{X_k-x}{|X_k-x|}\in\{\nu,\nu'\}\}.
\end{align*}
We can then use the Markov property to show that  $\mathbb{P}-a.s.$,
for any $n\geq 1$,
\begin{align*}
P_{0,\omega}((D^{\ell_0})^n<\infty,(D^{\ell_0})^{n+1}&=\infty,B_{\nu}\cup B_{\nu'})=\\
&\sum_{x\in\mathbb{Z}^d}P_{0,\omega}((D^{\ell_0})^{n}<\infty,X_{(D^{\ell_0})^n}=x)
P_{x,\omega}(D^{\ell_0}=\infty,B_{\nu}\cup B_{\nu'}).
\end{align*}
We finally obtain $P_0((D^{\ell_0})^n<\infty,\forall n \geq 1|B_{\nu}\cup B_{\nu'})=1$.
This establishes a contradiction with $\textit{(H')}$ and concludes the
proof of (\ref{positif}).
It is possible to choose $\alpha>0$ such that for every
 $i\in\llbracket 2,d \rrbracket$:
 \begin{equation*}
  \ell'_i\left(\alpha\right)\in\mathcal{V}\ \textrm{and}\  \ell'_{-i}\left(\alpha\right)\in\mathcal{V},
 \end{equation*}
in the notations of (\ref{eq:defli}) with $\ell=\ell_0$.
From (\ref{eq:cond}) and Proposition $1.2$ in \cite{SZ_LLN}, it is
clear that almost surely on $\{B_{\nu}\cup B_{\nu'}\}\cap\{D^{\ell_0}=\infty\}$,
\begin{equation*}
 \tau^{\ell_2'\left(\alpha\right)}\vee\dots\vee\tau^{\ell_d'\left(\alpha\right)}\vee\tau^{\ell'_{-2}\left(\alpha\right)}\vee\dots\vee\tau^{\ell'_{-d}\left(\alpha\right)}<\infty.
 \end{equation*}
We can then follow the end of the proof of Lemma \ref{lemme:cone} and we obtain:
 \begin{equation*}
\union_{\alpha>0}\{D_{\alpha}^{\ell_0}=\infty,B_{\nu}\cup B_{\nu'}\}\overset{P_0- a.s.}{=}\{D^{\ell_0}=\infty,B_{\nu}\cup B_{\nu'}\},
\end{equation*}
and hence we can fix $\alpha>0$ such that,
\begin{equation*}
P_0(D^{\ell_0}_{\alpha}=\infty,B_{\nu}\cup B_{\nu'})>0.
\end{equation*}
We will now show that on $A_{\ell_0}$ (this is much easier
than on $B_{\nu}\cup B_{\nu'}$), the walk admits almost surely a unique asymptotic
direction. We use the same renewal structure with cones as in the
proof of Theorem \ref{theorem1}. Following the proof of Proposition
$1.2$ in \cite{SZ_LLN}, we obtain
\begin{proposition}Under assumption \textit{(H')},
\begin{equation*}
P_0-a.s., \qquad A_{\ell_0}=\{K<\infty\}=\{\tau_1<\infty\}.
\end{equation*}
\end{proposition}
We also adopt the notation $Q_0$ to denote the probability
measure $P_0(\cdot|A_{\ell_0})$. The proof of Corollary $1.5$
in \cite{SZ_LLN} leads to, 
\begin{proposition}Under assumption \textit{(H')},\\
% \label{prop:renouvellement2}
 $ \left( \left(X_{ \tau_1 \wedge\cdot} \right),\tau_1
   \right),\left(\left(X_{ \left( \tau_1+ \cdot \right) \wedge
   \tau_2}-X_{\tau_1}\right), \tau_2-\tau_1
   \right),\dots, \left( \left( X_{ \left( \tau_k + \cdot \right)
   \wedge \tau_{k+1}}-X_{\tau_k}\right),\tau_{k+1}-\tau_k
   \right)$
 are independent variables under $Q_0$ and for $k\geq 1$,
 $\left(\left(X_{\left(\tau_k+\cdot\right)\wedge\tau_{k+1}}-X_{\tau_k}\right),\tau_{k+1}-\tau_k\right)$
 are distributed like
 $\left(\left(X_{\tau_1\wedge\cdot}\right),\tau_1\right)$ under
 $P_0\left(\cdot|D^{\ell_0}_{\alpha}=\infty\right)$.
 \end{proposition}

Using again Lemma $3.2.5$ p$265$ in \cite{ZEI}, we also have
\begin{lemma}Under assumption \textit{(H')},
\begin{equation*}
E_{Q_0}[X_{\tau_1}\cdot\ell|D^{\ell}_{\alpha}=\infty]=\frac{1}{P_0\left(D_{\alpha}^{\ell}=\infty\right)}.
\end{equation*}
\end{lemma}
We have now all the tools to follow the proof of Theorem
\ref{theorem1} except that $Q_0$ replaces $P_0$. We obtain the
existence of $\nu_0$ satisfying (\ref{but}).
\end{proof}

%%%%%%%%%%%%%%%%%%%%%%%%%%%%%%%%%%%%%%%%%%%%%%%%%%%%%%%%%%%%%%%%%%%%%%%%%%%%%%%%%%%%%%%%%%%%%%%%%
%%%%%%%%%%%%%%%%%%%%%%%%%%%%%%%%%%%%%%%%%%%%%%%%%%%%%%%%%%%%%%%%%%%%%%%%%%%%%%%%%%%%%%%%%%%%%%%%%

\textbf{Acknowledgments}: I wish to thank my Ph.D. supervisor
Francis Comets for his help and suggestions. I am also grateful to the
anonymous referee for his remarks and comments.

%%%%%%%%%%%%%%%%%%%%%%%%%%%%%%%%%%%%%%%%%%%%%%%%%%%%%%%%%%%%%%%%%%%%%%%%%%%%%%
%%%%%%%%%%%%%%%%%%%%%%%%% The bibliography %%%%%%%%%%%%%%%%%%%%%%%%%%%%%%%%%%%
%%%%%%%%%%%%%%%%%%%%%%%%%%%%%%%%%%%%%%%%%%%%%%%%%%%%%%%%%%%%%%%%%%%%%%%%%%%%%%

\bigskip

\bibliographystyle{plain}
\bibliography{article}
% \begin{thebibliography}{20}

% %%%%%%%%%% The immortal Feller %%%%%%%%%%%%%
% \bibitem{cf:Szn98} Sznitman, A.S.: \emph{Brownian Motion, Obstacles
%     and Random media}, Springer monographs in mathematics, Springer (1998).

% \end{thebibliography}

\end{document}